\documentclass[11pt]{amsart}

\usepackage[final]{nick}

\newcommand{\schwartz}[2]{\calS_{#1}[#2]}
\newcommand{\euclid}{\calE(3)}
\newcommand{\euclidp}{\calE^+(3)}
\newcommand{\sgauge}[2]{G_{#1}[#2]}

\newlength{\myVSpace}
\newcommand{\LoopR}{\Lambda_{(r,1)}\matSL{2}{\bbC}}
\newcommand{\barbar}[1]{{{#1}^\circ}}
\DeclareMathOperator*{\mult}{mult}

\title{Constant Mean Curvature $n$-noids with Platonic Symmetries}

\author{N. Schmitt}

\begin{document}



\maketitle




\typeout{== intro.tex =============================================}
\section*{Introduction}

We construct families of genus-zero
constant mean curvature (CMC) surfaces in $\bbR^3$
with finite symmetry groups.
\Note{Give a brief history of the problem, listing partial and numerical
results.}

The construction is via the extended Weierstrass representation
for constant mean curvature surfaces.
A description of this construction
(Dorfmeister-Pedit-Wu~\cite{Dorfmeister_Pedit_Wu_1998})
can be found in~\cite{Schmitt_Kilian_Kobayashi_Rossman_2006},
briefly in~\autoref{sec:weierstrass}, and elsewhere.
The construction for immersions into $\bbR^3$ is outlined as follows.
(1)~Write down a potential $\xi$ with values in the $\matsl{2}{}$-valued
loops on the $\lambda$-parametrized unit circle $\bbS^1$ and
certain meromorphic behavior in $\lambda$ on the unit disk.
(2)~Solve the equation $d\Phi=\Phi\xi$.
(3)~Let $F$ be the unitary factor in the Iwasawa factorization of $\Phi$.
Then $F$ is the extended frame of a local CMC immersion.
The immersion can be obtained from $F$ via a Sym formula.

The ends of the immersions are conjectured to be
asymptotic to half-Delaunay surfaces,
because their potentials are local perturbations of Delaunay
potentials at their ends.

Since the construction of the immersion involves
solving a differential equation which has monodromy,
the difficulty of the construction is to close the ends.
This amounts to showing that the monodromy representation
is pointwise simultaneously unitarizable
on the unit circle via conjugation by a dressing.

The unitarization problem for the immersions in this paper is solved
by reduction to the trinoid case, which
is solved~\cite{Schmitt_Kilian_Kobayashi_Rossman_2006}. The idea is to
impose symmetries on the potential
in such a way that the fundamental piece of the surface
has at most three points ends.
A modification of the trinoid technology for non-closed surfaces
then shows that under a suitable set of inequalities,
the monodromy for the fundamental piece is unitarizable.
The monodromies for the whole immersion are powers of the monodromies
of the fundamental piece, so that the monodromy group for
the whole immersion is unitarizable.

More precisely, the above program is implemented as follows.
We construct a rational map $u:\bbP^1_z\to\bbP^1_u$
together with a potential $\xi$ on $\bbP^1_z$
such that $\xi$ descends via $u$ to a potential
$\eta$ on $\bbP^1_u$ which has at most three poles
at the branch values of $u$.
This is accomplished by taking $u$ to be an invariant of finite
group of M\"{o}bius transformations acting on $\bbP^1_z$.
By trinoid technology,
the spherical triangle inequalities on the
logs of the eigenvalues of a set of generators of the monodromy group
imply that the monodromy for $\eta$ is unitarizable.
The gauge equivalence between $\eta$ and $\xi$ then implies that
the monodromy for $\xi$ is
unitarizable. It follows that $\xi$ induces a CMC immersion of a punctured
Riemann sphere.
\Note{The pictures in this draft are from an old draft.
The final paper will have a new set of pictures,
including the platonic guys with ends on the face centers
and surfaces in spaceforms.}
\Note{Decide whether to christen these surfaces ``astronoids,'' or stick
to ``$n$-noids.''}

Since the astronoid construction is based on the trinoid construction,
many of the tools have already been developed~\cite{Schmitt_Kilian_Kobayashi_Rossman_2006}.
The paper starts with lemmas on on
finite symmetry groups and invariants,
and on monodromy for potentials and their pullbacks.
\autoref{sec:unitarize}
solves the astronoid monodromy unitarizability problem,
and the astronoids are constructed in \autoref{thm:nnoid}.
The remainder of the paper is concerned with symmetry.
\autoref{thm:symmetry-closed} shows that
the symmetries of the potential induce symmetries
of the CMC immersion, under suitable conditions.
The symmetries of the astronoids are
computed in \autoref{thm:nnoid-symmetry}.

I thank Franz Pedit for many helpful discussions.

\typeout{== weierstrass.tex =============================================}
\section{The extended Weierstrass representation}
\label{sec:weierstrass}

\subsection{Notation}
We will use the following notations:

For $r\in(0,\,1]$, $\calC_r=\{\lambda\in\bbC\suchthat\abs{\lambda}=1\}$.
For $r\in(0,\,1]$, $\calD_r=\{\lambda\in\bbC\suchthat \abs{\lambda}<r\}$.
For $r\in(0,\,1)$, $\calA_r=\{\lambda\in\bbC\suchthat r<\abs{\lambda}<1/r\}$.

For $r\in(0,\,1]$,
$\LoopSL{r}$ is the group of analytic maps $\calC_r\to\matSL{2}{\bbC}$.

For $r\in(0,\,1)$,
$\LoopuSL{r}\subset\LoopSL{r}$ is the subgroup of $r$-unitary loops,
that is, analytic maps $\calC_r\to\matSL{2}{\bbC}$
which are the boundaries of maps $\calA_r\to\matSL{2}{\bbC}$
satisfying the
reality condition $X^\ast = X^{-1}$, where $X^\ast$ is defined by
\[
X^\ast(\lambda) = \transpose{\ol{X(1/\ol{\lambda})}}.
\]
$\LoopuSL{1}\subset\LoopSL{1}$ is the group of analytic maps
$\bbS^1\to\matSU{2}{}$.

For $r\in(0,\,1]$, $\LooppSL{r}$ is the group of analytic maps
$\calC\to\matSL{2}{\bbC}$ which are the boundaries of analytic
maps $\bbD_r\to\matSL{2}{\bbC}$, and such that $X(0)$ is upper-triangular
with real entries on its diagonal.

$\LoopR$ is the group of holomorphic maps from the
annulus $\{\lambda\in\bbC\suchthat r <\abs{\lambda}<1\}$
to $\matSL{2}{\bbC}$.

For $r\in(0,\,1]$,
the $r$-Iwasawa factorization~\cite{McIntosh_Nonlinearity_1994, Burstall_Pedit_1995} of $X\in\LoopSL{1}$ is the unique
factorization $X=FB$, where $F\in\LoopuSL{r}$ and $B\in\LooppSL{r}$.

For $D\subset\bbC$ a domain invariant under complex conjugation
and $X:D\to\matSL{2}{\bbC}$ a holomorphic map, we will use
the notation
$\barbar{X(\lambda)}=\ol{X{\ol{\lambda}}}$.

\subsection{The extended Weierstrass construction}

The extended Weierstrass representation for CMC surfaces
is as follows:

1. Let $\Sigma$ be a domain in $\bbC$. With $r\in(0,\,1]$, let $\xi$ be a
\emph{potential}, that is, a
$\Loopsl{r}$-valued
differntial form which is the boundary of a meromrophic differnential
on the unit $\lambda$ disk, with a pole only at $\lambda=0$, and only
in the upper-right entry.

2. Solve the ordinary differential equation $d\Phi=\Phi\xi$.

3. Let $F$ be the unitary factor in the $r$-Iwasawa factorization of $\Phi$.
Then $F$ is the extended frame of some CMC immersion.

4. Apply the Sym formula $f=F'F^{-1}$,
where the prime denotes differentiation with respect to $\theta$,
where $\lambda=e^{i\theta}\in\bbS^1$.

If $\calM$ is the monodromy group, of $F$,
then the immersion $f$ on the universal
cover descends to an immersion on $\Sigma$ if and only iff all $M_F\in\calM$
satisfy
\begin{equation}
\label{eq:closing}
M_F(1)=\pm\id,\quad M_F'(1)=0.
\end{equation}

\typeout{== algebra.tex =============================================}
\section{Finite symmetry groups and invariants}
\label{sec:algebra}

\subsection{\texorpdfstring{Finite subgroups of $\matSO{3}{}$}{Finite subgroups of SO3}}

We give here for reference the well-known classification
of the finite subgroups of $\matSO{3}{}$.

Let $\euclid$ denote the isometries of Euclidean 3-space,
and $\euclidp\subset\euclid$ the subgroup of orientation-preserving
isometries.

The finite subgroups of $\matPSU{2}{}\cong\matSO{3}{}$
and $\matPSL{2}{\bbC}$ are classified
as follows~\cite{Coxeter}:

1. Every finite nontrivial subgroup of $\matPSU{2}{}\cong\matSO{3}{}$
is one of the following groups:
\begin{itemize}
\item the cyclic group $\bbZ_n$ (order $n$),
\item the dihedral group $D_n$ (order $2n$),
\item the symmetry group $A_4$ of the tetrahedron (order $12$),
\item the symmetry group $S_4$ of the octahedron or hexahedron (order $24$),
\item the symmetry group $A_5$ of the icosahedron or dodecahedron (order $60$).
\end{itemize}

2. If two finite subgroups of $\matPSU{2}{}$ are isomorphic, then they
are conjugate via an element of $\matPSU{2}{}$.

3. Every finite subgroup of $\matPSL{2}{\bbC}$ is isomorphic
to a subgroup of $\matPSU{2}{}$ via conjugation by an element
of $\matPSL{2}{\bbC}$.

4. Every finite subgroups of $\euclid$ is isomorphic to a subgroup
of $\matO{3}{}$ via a conjugation by an element of $\euclid$.

5. Every finite subgroups of $\euclidp$ is isomorphic to a subgroup
of $\matSO{3}{}$ via a conjugation by an element of $\euclidp$.

Statement 3 can be proved by
averaging an inner product on $\bbC^2$ with respect to
the action of the finite subgroup $G$ to get an inner product
which is invariant under the action of $G$.

\subsection{The rational invariant}

For each finite subgroup $\subset\matPSL{2}{\bbC}$ of order $d$, there
is an essentially unique rational function $u$ of degree $d$
which is invariant under the action of $G$ and which generates the field
of such invariant functions.


Let $\bbC(z)$ denote the field of rational functions on $\bbP^1$
in the variable $z$.
For a group $G\subset\matPSL{2}{\bbC}$, let
$\bbC(z)^G\subset\bbC(z)$ denote the field of rational functions
on $\bbP^1$ in the variable $z$ which are invariant under the
action of $G$, in the
sense that for all $v\in\bbC(z)^G$ and all $g\in G$, $v\circ g = v$.


\begin{lemma}
\label{lem:invariant}
For any subgroup $G\subset\matPSL{2}{\bbC}$ with finite order $n$,
the field $\bbC(z)^G$ of rational maps
invariant under $G$
is generated by a single element with degree $n$.
Moreover, given one such generator $u$,
any other generator is of the form $\tau\circ u$ for some $\tau\in\matPSL{2}{\bbC}$.
\Note{Prove uniqueness.}
\end{lemma}

\begin{proof}
To construct a generator, let $p,\,q\in\bbP^1$ be distinct points,
and define $u\in\bbC(z)$ by
\[
u(z) = \frac{\prod_{g\in G}(z-g(p))}{\prod_{g\in G}(z-g(q))}.
\]
To show that $u\in\bbC(z)^G$, take $h\in G$. Then
\[
(u\circ h)(z) = \frac{\prod_{g\in G}(h(z)-g(p))}{\prod_{g\in G}((h(z)-g(q))}.
\]
The zero set (resp.~pole set) of $u\circ h$ is
the orbit of $p$ (resp.~$q$) under the action of $G$,
counted with multiplicities.
It follows that $u\circ h$ is a constant multiple of $u$.
At a fixed point $r$ of $h$, we have $(u\circ h)(r) = u(r)$,
so $u\circ h = u$.
This shows $u\in\bbC(z)^G$.

To show that $u$ generates $\bbC(u)^G$,
first note that the image of the orbit under $G$ of any $p\in\bbP^1$
under any element of $\bbC(u)^G$ is a single point.
It follows that the degree of every non-constant element of $\bbC(u)^G$
is a multiple of $n$.

Let $v\in\bbC(z)^G$.
Then $v$ takes the orbit of $p$ to $a$ and the orbit of $q$ to $b$.
Let $\rho\in\matPSL{2}{\bbC}$ satisfy $\rho(a)=0$, $\rho(b)=\infty$.
Then $v\circ h$ has zeros (resp.~poles) at the orbit of $p$ (resp.~$q$)
counted with multiplicities. Hence $\deg( (v\circ h)/u ) = \deg(v)-n$.
By induction, $v$ can be reduced to degree $>n$, implying degree $0$.
It follows by induction that $v$ is a rational function of $u$.
\end{proof}

Given a analytic function $f$ in a neighborhood of $p\in\bbC$
with series expansion $a_0 + a_j (z-p)^j + \Order((z-p)^{j+1})$,
$a_j \ne 0$,
we define $j$ as the \emph{multiplicity} $j = \mult_p f$ of $f$ at $p$.
\Note{This is one plus the branch order of $f$ at $p$?}
Note that $\mult_p f = 1 + \ord_p f'$.

\begin{lemma}
Let $G\subset\matPSL{2}{\bbC}$ be a finite subgroup of order $d$
and let $u$ be a generator of $\bbC(z)^G$ of degree $d$
as in \autoref{lem:invariant}.
Then
\begin{enumerate}
\item
The orbits of $G$ are the sets of the form $u^{-1}(\{x\})$, $x\in\bbP^1$.
\item
If $A\subset\bbP^1$ is the orbit of $a\in\bbP^1$ under $G$,
then $\mult_a u = d/\card A$.
\end{enumerate}
\end{lemma}



\begin{lemma}
\label{lem:orbit}
Let $G\in\matPSL{2}{\bbC}$ be a finite subgroup, and let
$u$ be a generator of $\bbC(z)^G$.
(1) $u$ has three branch values, except in the case $G=\bbZ_n$,
in which it has two.
(2) If $u(p)=u(q)$, then $\mult_p u = \mult_q u$.
(3)
The following table lists, for each finite subgroup group
$G\subset\matSL{2}{\bbC}$, the multiplicities of $u$
at the branch points of $u$ which map to distinct branch values.
\[
\begin{array}{c|l}
\text{Group} & \text{Multiplicities}\\
\hline
\bbZ_n & n,\,n \\
D_n    & 2,\,2,\,n \\
A_4    & 3,\,3,\,2\\
S_4    & 4,\,3,\,2\\
A_5    & 5,\,3,\,2
\end{array}
\]
\end{lemma}

\begin{proof}
For the case $G=\bbZ_n$, generated by $z\mapsto e^{2\pi i/n}z$, then
$u=z^n$ is a generator of $\bbC(z)^G$.
Then the branch points of $u$ are $0$ and $\infty$ with
respective branch values $0$ and $\infty$, and
$\mult_0 u = \mult_\infty u = n$.

Assume now that $G\ne\bbZ_n$.
Then $G$ is the orientation-preserving
symmetry group of a double right pyramid,
tetrahedron, octahedron or icosahedron.
Let $P$ be the centrally-projected polyhedron to the 2-sphere.
Let $B_1$ be the set of vertices, $B_2$ the set of face centers,
and $B_3$ the set of edge centers.
(For the dihedral group, let $B_1$ be the set of vertices other than
the apexes, $B_2$ the set of two apexes, and $B_3$ the edge
centers of the $n$-gon between the apexes.)
Then it is geometrically clear that $B_1$, $B_2$, and $B_3$ are
distinct orbits of $G$, with cardinalities
$(\card B_1,\,\card B_2,\,\card B_3)$ respectively
$(n,\,n,\,2)$,
$(4,\,4,\,6)$,
$(6,\,8,\,12)$,
$(12,\,20,\,30)$.

By \autoref{lem:orbit}, $u$ has the multiplicities listed in the
table above.

To show that $u$ has no other branch points we use the formula
\[
\sum_{p\in\bbP^1} (\mult_p u-1) = 2d-2,
\]
where $d$ is the order of $G$ and the degree of $u$.
For $k=1,\,2,\,3$, define $n_k = \mult_p u$, where $p\in B_k$.
Let $B=B_1\cup B_2\cup B_3$.
Then
\[
\sum_{p\in B}(\mult_p u-1)
 = \sum_1^3 \card B_k(n_k-1)
 = 3d-\sum_1^3 \card  B_k = 2d-2.
\]
Hence there cannot be any points in $\bbP^1\setminus B$
such that $\mult_p u>1$, so $u$ has no branch points in $\bbP^1\setminus B$.
\end{proof}

A generator $u$ of $\bbC(z)^G$ also has reflective symmetries.

\begin{lemma}
\label{lem:reflectu}
Let $G\subset\matPSL{2}{\bbC}$ be a finite subgroup of order $n$ which is not
cyclic, let $u$ be a generator of $\bbC(z)^G$ of degree $n$,
let $b_1,\,b_2,\,b_3$ be the branch values of $u$, and let
$B_k = u^{_1}(\{b_k\})$, $k=1,\,2,\,3$.
Let $\sigma$ be an inversion in a circle which preserves each of the
sets $B_1,\,B_2,\,B_3$. Let $\tau$ be the inversion in the circle
through $b_1,\,b_2,\,b_3$.
Then $\tau\circ u\circ \sigma = u$.
\end{lemma}

\begin{proof}
Let $f=\tau\circ u\circ \sigma$. Then $f$ is a rational function on
$\bbP^1$ with the same degree as $u$.
Then for each $k\in\{1,\,2,\,3\}$, $f^{-1}(\{b_k\}) = B_k$.
Moreover, $f$ has the same multiplicity as $u$ at each element of
$B_1\cup B_2\cup B_3$.\Note{Why?}
Since the cardinality of $B_1\cup B_2\cup B_3$ counted with multiplicity
is $3n$, and $u$ is determined by its value at $2n+2$ points,
then $f=u$.
\end{proof}

%
\begin{lemma}
\label{lem:u}
Let $G_0\subset\matPSL{2}{\bbC}$ be a finite subgroup.
Then there exists a subgroup $G\subset\matPSL{2}{\bbC}$
conjugate to $G_0$ by an element of $\matPSL{2}{\bbC}$,
and a generator $u$ of $\bbC(z)^G$,
such that the branch values of $u$ are $0,\,1,\,\infty$
(or $0,\,\infty$ in the case $G=\bbZ^n$),
and $u$ satisfies $\ol{u(\ol{z})}=u(z)$.
\end{lemma}

\begin{proof}
In the case $G=\bbZ_n$, $G_0$ is conjugate to the group $G$ generated
by $\tau\in\matPSU{2}{\bbC}$ defined by $\tau(z)=e^{2\pi i/n} z$.
Then $u(z)=z^n$ is the a generator satisfying the conditions.

Hence assume $G\ne\bbZ_n$.
Let $u_0$ be a generator of $\bbC(z)^G$ and let
$b_k$ ($k=1,\,2,\,3$), $\sigma$ and $\tau$ be
as in \autoref{lem:reflectu}, so that
by that lemma, $u=\tau\circ u\circ\sigma$.
By a conjugation, we may assume that $\tau$ is the map $\tau(z)=\ol{z}$.
By a M\"obius transformation, we may assume that
$(b_1,\,b_2,\,b_3) = (0,\,1,\,\infty)$.
Then $u$ satisfies the required conditions.
\end{proof}

\section{The Schwartz gauge}

\typeout{== potential.tex =============================================}
\subsection{The potential}
\label{sec:potential}

The construction of CMC surfaces by the extended Weierstrass representation
begins with a meromorphic $\matsl{2}{\bbC}$-valued
potential on the underlying Riemann surface.
Let $z$ be the standard conformal coordinate on the Riemann sphere
$\bbP^1$.
We will use a potential of the form
\begin{equation}
\label{eq:potential}
\xi = \begin{pmatrix} 0 & \lambda^{-1}dz \\
  (1-\lambda)^2 Q/dz & 0
\end{pmatrix}
\end{equation}
where $Q$ is the Hopf differential.

This choice of potential is motivated by the following considerations:

1. We would like the resulting immersion to have asymptotically Delaunay
ends. Hence all poles of $Q$ should be double poles, with real
quadratic residues.

2. Since we aim for Delaunay ends, we want the eigenvalues of the
monodromy associated to $\xi$ to be the same as those for a Delaunay surface.
The form of $p$ insures this.

3. The form of $p$ insures that the closing conditions~\eqref{eq:closing}
hold for the monodromy of the holomorphic frame.

4. Moreover, since ambient symmetries of the immersion induce symmetries
in the Hopf differential, we want to choose $Q$ to be invariant under
whatever symmetry group we want to impose on the surface.

\typeout{== schwartz.tex =============================================}
In this section we define and determine some properties of
the \emph{Schwartz gauge}, a gauge $g$ of an off-diagonal potential
on a Riemann sphere whose upper-right term $\lambda^{-1}dz$
is left invariant under the gauge action of $g$.

\subsection{The Schwartzian derivative}


For a function $u=u(z):\bbP^1\to\bbP^1$,
the Schwartzian derivative $\schwartz{z}{u}$ of $u$ with respect to $z$
is the quadratic differential on $\bbP^1$ defined by
\begin{equation}
\label{eq:schwartz}
\frac{\schwartz{z}{u}}{dz^2}
 = \biggl(\frac{u''}{u'}\biggr)' - \frac{1}{2}\biggl(\frac{u''}{u'}\biggr)^2,
\end{equation}
where the prime notation denotes differentiation with respect to $z$.

\begin{lemma}
\label{lem:schwartz}
\mbox{}

\begin{enumerate}
\item
For all $\tau\in\matPSL{2}{\bbC}$, $\schwartz{z}{\tau}=0$.
\item
For $u\in\bbC(w)$ and $v\in\bbC(z)$,
the ``Schwartz derivative chain rule'' holds:
\begin{equation}
\label{eq:schwartz-chain-rule}
\schwartz{w}{v\circ u} = u^\ast(\schwartz{z}{v}) + \schwartz{w}{u}.
\end{equation}
\item
Let $v\in\bbC(z)$ and $\tau\in\matPSL{2}{\bbC}$.
If $\tau^\ast v=v$,
then $\tau^\ast(\schwartz{z}{v}) = \schwartz{z}{v}$.
\item
Let $G\subset\matPSL{2}{\bbC}$ be a finite group and let
$u$ be a generator of $\bbC(z)^G$.
Then there exists a quadratic differential $\alpha$ on $\bbP^1_u$ such
that $v^\ast\alpha = \schwartz{z}{u}$.
\end{enumerate}
\end{lemma}

\begin{proof}
Statements (i)--(iii) are computations.

To prove (iii),
\Note{Clean up $du$ and $u^\ast du$ confusion.}
note that by (ii), for all $\tau\in G$,
$\tau^\ast\schwartz{z}{u}=\schwartz{z}{u}$.
Consider $f=\schwartz{z}{u}/du^2$.
Then for all $\tau\in G$, we have
\[
\tau^\ast f
 = \tau^\ast \schwartz{z}{u} / \tau^\ast du^2
 = \schwartz{z}{u}/du^2 = f.
\]
Hence $f\in\bbC(z)^G$.
Since $u$ generates $\bbC(z)^G$, then
there exists $g$ such that $u^\ast g = f$.
Let $\alpha = g du^2$.
Then
\[
u^\ast \alpha = (u^\ast g)(u^\ast du^2) = f du^2 = \schwartz{z}{u}.
\]
\end{proof}

\typeout{== gauge.tex =============================================}
\subsection{The Schwartz gauge}


\begin{definition}[Schwartz gauge]
Let $\Sigma_1,\,\Sigma_2$ be Riemann surfaces and suppose
there exists a global coordinate $z$ on $\Sigma_1$.
Let $u:\Sigma_1\to\Sigma_2$ be a map.
Let $B\subset\Sigma_1$ be the branch points of $u$.
Then we define the \emph{Schwartz gauge}
$\sgauge{z}{u}:(\Sigma_1\setminus B)\times\bbC_\lambda\to\matSL{2}{\bbC}$ by
\[
\sgauge{z}{u} =
  \begin{pmatrix}v & 0 \\ -\lambda v' & v^{-1}\end{pmatrix},\quad
v={(u')}^{-1/2},
\]
where the prime denotes differentiation with respect to $z$.

Note that while $\sgauge{z}{u}$ is only defined up to sign,
gauge action by $\sgauge{z}{u}$ is well-defined.
\end{definition}

\begin{lemma}
\label{lem:gauge1}
Let $z$ be the standard coordinate on $\bbP^1\setminus\{\infty\}$.
Let $Q$ be a quadratic differential on $\Sigma$.
Let $u:\bbP^1\to\bbP^1$ be a rational function
let $B$ be the set of branch points of $u$,
let $\Sigma=\bbP^1\setminus B$,
and let $\widetilde{\Sigma}$ be a double cover of $\Sigma$.
Let $\xi$ be a potential of the form~\eqref{eq:potential}.
Then
\begin{equation}
\label{eq:schwartz-gauge}
\gauge{\xi}{\sgauge{z}{u}} = 
\begin{pmatrix}0 
& \lambda^{-1}du \\
 \bigl((1-\lambda)^2Q + \half \lambda\schwartz{z}{u}\bigr)/du
 & 0\end{pmatrix}.
\end{equation}
\end{lemma}

\begin{proof}
The proof is a direct calculation.
Let the prime denote differentiation with respect to $z$,
so $u'(z)dz=du$. Then by~\eqref{eq:schwartz},
\[
\schwartz{z}{u}
=
\left(\frac{u'''}{u'} - \frac{3}{2}\biggl(\frac{u''}{u'}\biggr)^2\right)dz^2
=
\left(\frac{u'''}{u'} - \frac{3}{2}\left(\frac{u''}{u'}\right)^2\right)\frac{du^2}{ {(u')}^2 }.
\]
Then with $g=\sgauge{z}{u}$,
\begin{align*}
\gauge{\xi}{g} &= g^{-1}\xi g + g^{-1}dg\\
&= \begin{pmatrix} \frac{u''}{2 u'} & u' \\
(1-\lambda)^2\frac{Q}{u'dz^2} - \frac{ {(u'')}^2}{4{(u')}^3} &
 -\frac{u'' dz}{2u'}\end{pmatrix}dz
+
\begin{pmatrix}
-\frac{u''}{2u'} & 0 \\
 -\frac{ {(u'')}^2}{2 {(u')}^3} + \frac{u'''}{2{(u')}^2} & \frac{u''}{2u'}
\end{pmatrix}dz\\
&=
\begin{pmatrix}
0 & u' \\
(1-\lambda)^2\frac{Q}{u'dz^2} - \frac{3 {(u'')}^2}{4 {(u')}^3 } + \frac{u'''}{2{(u')}^2} & 0
\end{pmatrix}dz\\
&=
\begin{pmatrix}
0 & du \\
(1-\lambda)^2\frac{Q}{du} +\frac{1}{2{(u')}^{2}}\left(\frac{u'''}{u'}-\frac{3}{2}\left(\frac{u''}{u'}\right)^2\right)du & 0
\end{pmatrix}\\
&=
\begin{pmatrix}
0 & du \\
\frac{1}{du}\bigl( (1-\lambda)^2Q + \half\schwartz{z}{u}\bigr) & 0
\end{pmatrix}.
\end{align*}
\end{proof}


The next lemma gives the formula for the Schwartz gauge
for a composition of maps.

\begin{lemma}
\label{lem:gauge-composition}
Let $\Sigma_1,\,\Sigma_2,\,\Sigma_3$ be Riemann surfaces
and suppose there exist global coordinates $z,\,w$ on
$\Sigma_1$ and $\Sigma_2$ respectively.
Let $u:\Sigma_1\to\Sigma_2$ and $v:\Sigma_2\to\Sigma_3$ be maps.
Then for some choice of sign,
\begin{equation}
\label{eq:sgauge}
\sgauge{z}{v\circ u} = \pm \sgauge{z}{u} \cdot u^\ast \sgauge{w}{v}.
\end{equation}
\end{lemma}

\begin{proof}
Note: equality~\eqref{eq:sgauge} is valid only up to sign
because of the sign ambiguity of its terms.

The proof is by direct calculation.
Let
\[
A = {\left(\frac{du}{dz}\right)}^{-1/2}
\quad\text{and}\quad
B = {\left(\frac{dv}{dw}\right)}^{-1/2}.
\]
Then
\[
\sgauge{z}{u} =
 \begin{pmatrix} A & 0 \\ -\lambda \frac{dA}{dz} & A^{-1}\end{pmatrix}
\quad\text{and}\quad
\sgauge{w}{v} =
 \begin{pmatrix} B & 0 \\ -\lambda \frac{dB}{dw} & B^{-1}\end{pmatrix}\\
\]
Hence
\begin{equation}
\label{eq:sgauge5}
\sgauge{z}{u} \cdot u^\ast\sgauge{w}{v} =
 \begin{pmatrix}
  A u^\ast B & 0 \\
  -\lambda\left(u^\ast B \frac{dA}{dz} + A^{-1}u^\ast \frac{dB}{dw}\right) & A^{-1}u^\ast B^{-1}.
  \end{pmatrix}
\end{equation}
On the other hand,
\begin{equation}
\label{eq:sgauge3}
\sgauge{z}{v\circ u}
 =
 \begin{pmatrix}
  A u^\ast B & 0 \\
  -\lambda \frac{d}{dz}\left(A u^\ast B\right) & A^{-1}u^\ast B^{-1}
 \end{pmatrix}.
\end{equation}
Looking at the lower-left entry of~\eqref{eq:sgauge3},
\begin{equation}
\label{eq:sgauge4}
\frac{d}{dz}\left(A u^\ast B\right)
 = u^\ast B \frac{dA}{dz} + A \frac{d}{dz}u^\ast B.
\end{equation}
The second term of the right hand side of~\eqref{eq:sgauge4} is
\[
A \frac{d}{dz}u^\ast B
 = A \cdot u^\ast \frac{dB}{dw} \cdot \frac{du}{dz} = A^{-1} u^\ast\frac{dB}{dw},
\]
so~\eqref{eq:sgauge4} becomes
\[
\frac{d}{dz}\left(A u^\ast B\right) = 
u^\ast B \frac{dA}{dz} + A^{-1} u^\ast \frac{dB}{dw}.
\]
Hence the lower-right entries of~\eqref{eq:sgauge5} and~\eqref{eq:sgauge3}
are equal, proving~\eqref{eq:sgauge}.
\end{proof}

The next lemma shows how invariance of the Hopf differential
under the action of a conformal or anticonformal M\"obius transformation
implies gauge invariance of a potential of the form~\eqref{eq:potential}.

\begin{lemma}
\label{lem:gauge2}
Let $\xi$ be a potential of the form~\eqref{eq:potential}
with Hopf differential $Q$.

\begin{enumerate}
\item
Suppose $\tau^\ast Q = Q$ for some M\"obius transformation
$\tau\in\matPSL{2}{\bbC}$.
Then $\tau^\ast\xi = \gauge{\xi}{\sgauge{z}{\tau}}$.
\item
Suppose $\ol{\tau^\ast Q} = Q$
for some anticonformal M\"obius transformation $\tau$.
Then $\tau^\ast\barbar{\xi} = \gauge{\xi}{\sgauge{z}{\ol\tau}}$.
\end{enumerate}
\end{lemma}

\begin{proof}
Part (i) follows by \autoref{lem:gauge1} and the fact
that $\schwartz{z}{\tau}=0$.
Part (ii) follows similarly.
\end{proof}

\typeout{== monodromy.tex =============================================}
\section{Monodromy}

Before constructing the symmetric $n$-noids, we will need
some facts about monodromy.

\subsection{Monodromy and gauges}

\begin{lemma}
\label{lem:monodromy-conjugation}
Let $G$ be a Lie group.
Let $\Sigma$ be a Riemann surface.
Let $\xi$ be a potential on $\Sigma$ with values in the lie algebra
for $G$.
Let $\rho:\widetilde{\Sigma}\to\Sigma$ be the universal cover of $\Sigma$.
Let $\Delta$ be the group of deck transformations on $\widetilde{\Sigma}$.
Let $\Phi:\widetilde{\Sigma}\to G$ be a solution
to the equation $d\Phi = \Phi\xi$.
Let $\tau:\Sigma\to\Sigma$ be a conformal map,
and let $\tilde{\tau}:\widetilde\Sigma\to\widetilde\Sigma$
be a lift of $\tau$ to $\widetilde\Sigma$,
that is, $\tilde\tau$ is some map satisfying
$\tau\circ\rho = \rho\circ\tilde\tau$.
Let $\delta_1,\,\delta_2\in\Delta$ such that
$\delta_2\circ\tilde{\tau} = \tilde{\tau}\circ\delta_1$.
Then
\begin{enumerate}
\item
\[
M(\delta_1,\,\tilde{\tau}^\ast\Phi) = M(\delta_2,\,\Phi).
\]
\item
If for some $h,\,g\in G$ we have
$\tilde{\tau}^\ast\Psi = h\Phi g$,
and $g$ satisfies $\delta_1^\ast g = \pm g$, then
\[
M(\delta_2,\,\Psi) = \pm \Ad_h M(\delta_1,\,\Phi).
\]
\end{enumerate}
\end{lemma}

\begin{proof}
To prove (i), first note that
\[
\delta_1\tilde{\tau}^\ast = (\tilde{\tau}\circ\delta_1)^\ast
=(\delta_2\circ\tilde{\tau})^\ast = \tilde{\tau}^\ast\delta_2^\ast.
\]
Then
\begin{align*}
M(\delta_1,\,\tilde{\tau}^\ast\Phi)
&=
\delta_1^\ast\tilde{\tau}^\ast\Phi\cdot{(\tilde{\tau}^\ast\Phi)}^{-1}
=
\tilde{\tau}^\ast\delta_2^\ast\Phi\cdot{(\tilde{\tau}^\ast\Phi)}^{-1}\\
&=
\tilde{\tau}^\ast\left(\delta_2^\ast\Phi\cdot\Phi^{-1}\right)
=
\tilde{\tau}^\ast\left(M(\delta_2,\,\Phi)\right) = M(\delta_2,\,\Phi).
\end{align*}

To prove (ii), note that by part (i) we have
$M(\delta_2,\,\Psi)=M(\delta_1,\,\tilde{\tau}^\ast\Psi)$.
Then
\begin{align*}
M(\delta_2,\,\Psi)
&=
M(\delta_1,\,\tilde{\tau}^\ast\Psi)
=M(\delta_1,\,h\Phi g)\\
&=
\delta_1^\ast(h\Phi g) \cdot (h\Phi g)^{-1}
=h \cdot \delta_1^\ast\Phi ( \delta_1^\ast g \cdot g^{-1}) \Phi^{-1} h^{-1}\\
&=
\pm h ( \delta_1^\ast\Phi \cdot \Phi^{-1}) h^{-1}
=\pm h M(\delta_1,\,\Phi) h^{-1}
\end{align*}
\end{proof}

\subsection{Monodromy and covers}

In the following theorem we will use the following setup.

\begin{notation}
\label{not:cover}
Let $u:\bbP^1\to\bbP^1$ be a rational map.
Let $B$ be the set of branch points of $u$.
Define the punctured Riemann spheres
$\Sigma_1 = \bbP^1\setminus B$ and $\Sigma_2 = \bbP^1\setminus u(B)$.

Let $\xi$ be a potential on $\Sigma_1$ and $\eta$ a potential
on $\Sigma_2$, such that
\[
u^\ast\eta = \gauge{\xi}{\sgauge{z}{u}}.
\]

Let $\rho_1:\widetilde\Sigma_1\to\Sigma_1$
and $\rho_2:\widetilde\Sigma_2\to\Sigma_2$
be respective universal covers.
Let $\Delta_1$ and $\Delta_2$ be the respective deck transformation groups.

Let $\Phi:\Sigma_1\to\LoopR$ be a solution to the equation
$d\Phi = \Phi\xi$.
Let $\Psi:\Sigma_2\to\LoopR$ be the solution to the equation
$d\Psi = \Psi\eta$ such that
\[
u^\ast\Psi = \Phi \sgauge{z}{u}.
\]
\end{notation}

\begin{lemma}
\label{lem:monodromy-descent}
Consider \autoref{not:cover}.
\Note{Probably change how this is written.}
Let $\tau:\bbP^1\to\bbP^1$ be a M\"{o}bius transformation
such that $u\circ\tau=u$.
Then $u$ has finite order $n$.
By \autoref{lem:gauge1}
\begin{equation}
\label{eq:descend2}
\tau^\ast\xi= \gauge{\xi}{\sgauge{z}{\tau}}.
\end{equation}
Let $\tilde\tau:\Sigma_1\to\Sigma_1$ such that
$\tau\circ u = u\circ\tilde\tau$.
Hence
\begin{equation}
\label{eq:descend3}
\tilde\tau^\ast\Phi = h \Phi \sgauge{z}{\tau}
\end{equation}
for some $h\in\LoopR{}$.
\begin{enumerate}
\item
Let $\delta_2\in\Delta_2$ be a deck transformation such that
$\delta_2\circ u = u\circ\tau$.
Then
$M(\delta_2,\,\Psi) = h$.
\item
$\tilde\tau^n\in\Delta_1$ and
$M(\tilde\tau^n,\,\Phi) = h^n$.
\end{enumerate}
\end{lemma}

\begin{proof}
To prove (i), write $g=\sgauge{z}{u}$ and $p=\sgauge{z}{\tau}$. Then
\begin{equation}
\label{eq:descend4}
\begin{split}
M(\delta_2,\,\Psi)
 &= \delta_2^\ast\Psi \cdot \Psi^{-1}\\
 &= u^\ast\delta_2^\ast\Psi \cdot (u^\ast\Psi)^{-1}
 = \tau^\ast u^\ast\Psi \cdot (u^\ast\Psi)^{-1}\\
 &= \tau^\ast(\Phi g) \cdot (\Phi g)^{-1}
 = \tau^\ast\Phi \cdot \tau^\ast g \cdot g^{-1}\Phi^{-1}\\
 &= h \Phi \left(p \cdot \tau^\ast g \cdot g^{-1}\right)\Phi^{-1}.
\end{split}
\end{equation}
By \autoref{lem:gauge-composition},
\[
p \cdot \tau^\ast g
 =
\sgauge{z}{\tau} \cdot \tau^\ast \sgauge{w}{u}
= \sgauge{z}{u\circ\tau} = \sgauge{z}{u} = g.
\]
Hence
\[
p \cdot \tau^\ast g \cdot g^{-1} = \id,
\]
so the right-hand side of~\eqref{eq:descend4} is $h$,
proving (i).

To prove (ii),
write $g=\sgauge{z}{\tau}$ and $\delta_1 = {(\tilde\tau^n)}^\ast$.
Apply $\tilde\tau$ to~\eqref{eq:descend3} $n-1$ times to obtain
\[
\delta_1^\ast\Phi
=
{(\tilde\tau^n)}^\ast\Phi
 =
 h^n \Phi
 \left(g \cdot \tilde\tau^\ast g \cdot \dots \cdot {(\tilde\tau^{n-1})}^\ast g\right).
\]
By \autoref{lem:gauge-composition},
\[
 g \cdot \tilde\tau^\ast g \cdot \dots \cdot {(\tilde\tau^{n-1})}^\ast g
 = \sgauge{z}{\tau^n}
 = \id,
\]
so
$\delta_1^\ast\Phi = h^n \Phi$.
By definition $M(\delta_1,\,\Phi) = \delta_1^\ast\Phi \cdot \Phi^{-1}$,
so $M(\delta_1,\,\Phi)=h^n$, proving (ii).
\end{proof}

\typeout{== unitarize.tex =============================================}
\section{\texorpdfstring{$n$-noids}{n-noids}}
\label{sec:unitarize}

We are now prepared to construct the $n$-noids, in \autoref{thm:nnoid}.
The computation of their symmetry groups (\autoref{thm:nnoid-symmetry}
is deferred to the following \autoref{sec:symmetry} and \autoref{sec:main}.

\subsection{%
  \texorpdfstring{Unitarizing three matrices whose product is $\id$}%
  {Unitarizing three matrices whose product is I}}

The closing of the surface is acheived by unitarizing the monodromy.
To do this, we show that the monodromy group is a subgroup of
a monodromy group generated by three loops whose product is one.
This unitarization problem is solved as follows.

\begin{definition}
The \emph{spherical triangle inequalities} on
three numbers $\nu_1,\,\nu_2,\,\nu_3\in\bbR^3$ are
\begin{equation}
\label{eq:striangle}
\abs{\nu_1} + \abs{\nu_2} + \abs{\nu_3} \le 1,\quad
\abs{\nu_1} \le  \abs{\nu_2} + \abs{\nu_3},\quad
\abs{\nu_2} \le  \abs{\nu_1} + \abs{\nu_3},\quad
\abs{\nu_3} \le  \abs{\nu_1} + \abs{\nu_2}.
\end{equation}
The \emph{strict spherical triangle inequalities} are
the above inequalities with ``$\le$'' replaced by ``$<$.''
\end{definition}

The strict spherical triangle equalities provide necessary and sufficient
conditions for simultaneous unitarizability of three matrices
whose product is one~\cite{Schmitt_Kilian_Kobayashi_Rossman_2006}:

\begin{lemma}
\label{lem:trinoid}
Let $M_1,\,M_2,\,M_3\in\matSL{2}{\bbC}$ satisfy $M_1M_2M_3=\id$.
For $k=1,\,2,\,3$, let $\exp(\pm 2\pi i\nu_k)$ be
the eignevalues of $M_k$ with $\nu_k\in[-1/2,\,1/2]$. Then
$M_1,\,M_2,\,M_3$ are irreducible and simultaneously unitarizable
if and only if the strict spherical inequalities hold.
\end{lemma}

\subsection{Eigenvalues of the monodromy}
\label{sec:eigenvalues}

We compute the eigenvalues of a set of generators
of the monodromy group downstairs.

\begin{notation}
\label{not:nnoid}
Let $G$ be a subgroup of $\matPSU{2}\cong\matSO{3}{}$.
Let $u:\bbP^1_z\to\bbP^1_u$ be a generator of $\bbC(z)^G$
as in \autoref{lem:u},
chosen with branch values $0$, $1$, $\infty$
and reflective symmetry $\ol{u(\ol{z})}=u(z)$.
Let $n_0,\,n_1,\,n_\infty$ be the orders of $u$ at respective preimages
of $0,\,1,\,\infty$ under $u$, as listed at the end of
\autoref{sec:algebra}.
(In the case $G=\bbZ_n$, choose $u$ with branch values $0$, $\infty$,
and take $(n_0,\,n_1,\,n_\infty)=(n,\,1,\,n)$.)
For $w_0,\,w_1,\,w_\infty\in\bbR$,
define the Hopf differential on $\bbP^1_u$
\[
Q = \frac{a_0 + a_1 u + a_2 u^2}{u^2(u-1)^2}du^2,
\]
where the coefficients $a_0,\,a_1,\,a_2\in\bbR$ of the numerator
are uniquely determined by the condition
that the respective quadratic residues of $Q$ at $0,\,1,\,\infty$
are $w_0/16,\,w_1/16,\,w_\infty/16$.
\Note{Check whether this should be $w_k/(16 n_k)$.}
On $\bbP^1_z$ define the potential
\[
\xi =
\begin{pmatrix} 0 & \lambda^{-1}dz \\
 (1-\lambda)^2 (u^\ast Q)/dz & 0\end{pmatrix}.
\]

By \autoref{lem:schwartz}(iv),
there exists a quadratic differential $\alpha$ on $\bbP^1_u$
such that $u^\ast\alpha = \schwartz{z}{u}$.
Let $\eta$ be the potential on $\bbP^1_u$ defined by
\[
\quad
\eta =
\begin{pmatrix} 0 & \lambda^{-1}du \\
 ((1-\lambda)^2 Q + \half\lambda\alpha)/du & 0\end{pmatrix}.
\]
By \autoref{lem:gauge1}, $\xi$ and $\eta$ are related by the gauge
equation
\begin{equation*}
u^\ast\eta = \gauge{\xi}{\sgauge{z}{u}}.
\end{equation*}
\end{notation}

We now compute the eigenvalues of a set of generators of the
monodromy group for $\eta$ on $\bbP^1_u$.
The quadratic differential $\alpha$ has poles only at the
branch values of $u$, and since $u$ is meromorphic, these poles
are double poles with quadratic residues $n^{-2}-1$.
Hence the set of poles of $\eta$ is a subset of $\{0,\,1,\,\infty\}$.

\begin{lemma}
\label{lem:eigenvalues}
Let $u$, $\xi$ and $\eta$ be as in \autoref{not:nnoid}.
Let $M_0$, $M_1$, $M_\infty$ be the monodromies for $\eta$ around
$0$, $1$, $\infty$ respectively.
Then
$M_0$, $M_1$, $M_\infty$ generate the monodromy group for $\eta$,
and for $k=0,\,1,\,\infty$,
the eigenvalues of $M_k$ are $\exp(\pm 2\pi i \mu_k(w_k,\,\lambda))$, where
\begin{equation}
\label{eq:eigenvalue}
\mu_k(w,\,\lambda)
 = \half - \tfrac{1}{2n_k}\sqrt{1+ \lambda^{-1}(1-\lambda)^2 w/4},
\quad k=0,\,1,\,\infty.
\end{equation}
\end{lemma}

We now show that the monodromy group
generated by $M_0$, $M_1$, $M_\infty$
on the thrice-punctured sphere is pointwise unitarizable
on $\bbS^1$ except possibly at a finite subset.

For all the finite groups except the cyclic group,
for any weights near zero, the surface is balanced and
the monodromy is unitarizable.
In the case of the cyclic group, there is a balancing constraint:
the sum of the $n$ weights must be large enough to balance
the remaining two ends.

\begin{lemma}
\label{lem:unitarization}
Let $G\subset\matPSL{2}{}$ be a finite subgroup
and let $u$ be a generator of $\bbC(z)^G$.
Let $n_0,\,n_1,\,n_\infty$ be the respective orders of $u$ at the
preimages of its three branch values
(in the case $G=\bbZ_n$, let $(n_0,\,n_1,\,n_\infty) = (n,1,n)$).
Then there exists an open set $W\subset\bbR^3$ such that for all
$(w_1,\,w_2,\,w_3)\in W$,
the three functions $\mu_k(w_k,\,\lambda)$, $k=1,2,3$,
defined by~\eqref{eq:eigenvalue},  satisfy
the spherical triangles inequalites~\eqref{eq:striangle}
for all $\lambda\in\bbS^1$
and the strict spherical triangle inequalities
for all $\lambda\in\bbS^1\setminus\{1\}$.
\end{lemma}

\begin{proof}
We first prove the lemma in the case that $G\ne\bbZ_n$.
Define the four functions $f,\,f_0,\,f_1,\,f_\infty$ of $w=(w_1,\,w_2,\,w_3)$
and $\lambda$ by
\begin{align*}
f(w,\,\lambda) &= 1-\abs{\mu_0(w_1,\,\lambda)}-\abs{\mu_1(w_2,\,\lambda)}-\abs{\mu_\infty(w_3,\,\lambda)},\\
f_0(w,\,\lambda) &= -\abs{\mu_0(w_1,\,\lambda)} + \abs{\mu_1(w_2,\,\lambda)} + \abs{\mu_\infty(w_3,\,\lambda)}\\
f_1(w,\,\lambda) &= +\abs{\mu_0(w_1,\,\lambda)} - \abs{\mu_1(w_2,\,\lambda)} + \abs{\mu_\infty(w_3,\,\lambda)}\\
f_\infty(w,\,\lambda) &= +\abs{\mu_0(w_1,\,\lambda)} + \abs{\mu_1(w_2,\,\lambda)} - \abs{\mu_\infty(w_3,\,\lambda)}.
\end{align*}
These functions are defined so that
the strict (resp.~non-strict) triangle inequalities hold
for the $\mu_k(w_k,\,\lambda)$ if and only if
$f(w,\,\lambda)$ and the $f_k(w,\,\lambda)$ are positive
(resp.~non-negative).

It can be checked that for all $\lambda\in\bbS^1$,
$f(0,\,\lambda)>0$ and $f_k(0,\,\lambda)>0$, $k=1,\,2,\,3$.
Hence the three functions $\mu_k(0,\,\lambda)$, $k=1,\,2,\,3$ satisfy
the strict spherical triangles inequalites for all $\lambda\in\bbS^1$.
By the continuity of the $\mu_k$, $f$, and the $f_k$ ($k=1,\,2,\,3$),
there exists an open set $W\subset\bbR^3$ containing $0$ such that
$f$ and the $f_k$ are strictly positive
for all $w\in W$ and all $\lambda\in\bbS^1$.
This proves the lemma in the case $G\ne \bbZ_n$.

In the case $G=\bbZ_n$,
fix $v_0\in(0,\,4(1-1/n)/n)$ and let $w_0=(0,\,0,\,v_0)$.
Then the following can be checked:
(1) for all $\lambda\in\bbS^1$,
$f(w_0,\,\lambda) > 0$ and $f_\infty(w_0,\,\lambda) > 0$.
(2) For all $\lambda\in\bbS^1\setminus\{1\}$,
$f_0(w_0,\,\lambda) > 0$ and $f_1(w_0,\,\lambda) > 0$.
(3) For all $w\in\bbR^3$, $f_0(w,\,1)=0$ and $f_1(w,\,1)= 0$.
By the continuity of the $\mu_k$, $f$, and the $f_k$ ($k=1,\,2,\,3$),
there exists an open set $W\subset\bbR^3$ containing $w_0$ such that
for all $w\in W$,
$f$ and the $f_k$ are positive
for all $\lambda\in\bbS^1\setminus\{1\}$
and are non-negative at $\lambda=1$.
This proves the lemma in the case $G=\bbZ_n$.
\end{proof}

\begin{remark}
From the explicitly given functions $\mu_k$
in \autoref{lem:unitarization}
it is possible to explictly calculate the set $(w_1,\,w_2,\,w_3)$
for which the $\mu_k$ satisfy the
spherical triangle inequalities on $\bbS^1$.
For a similar calculation for trinoids (the case $n_0=n_1=n_\infty=1$),
see~\cite{Schmitt_Kilian_Kobayashi_Rossman_2006}.
\Note{I am currently working on this computation. If this
stronger result can be written reasonably briefly,
it will replace \autoref{lem:unitarization} and this
remark will disappear. The techniques for the trinoid computation
(the $n$'s all $1$)
do not trivially transfer to the present case.
The result should be: if the ``necksizes'' for the three downstairs monodromy
satisfy the spherical triangle inequalities, then they are
simultaneously unitarizable, and hence so is the upstairs monodromy.
The issue is showing that the spherical triangle inequalities
ont the three ``necksizes'' imply them for the logs of the monodromy
eignevalues on all of $\bbS^1$.}
\end{remark}


We now have the tools to construct CMC immersions of the punctured
sphere. The computation of their symmetries
is deferred to \autoref{thm:nnoid-symmetry} below.

\begin{theorem}
\label{thm:nnoid}
Let $G\subset\matSO{3}{}$ be a finite subgroup.
Then there exists an open set $W\subset\bbR^3$ such that
for each choice $(w_0,\,w_1,\,w_\infty)\in W$,
there exists an initial value $C\in\LoopR$ such that
the CMC immersion induced by
\[
d\Phi = \Phi\xi,\quad \Phi(b)=C
\]
lives on the punctured sphere.
\end{theorem}

\begin{proof}
By \autoref{lem:unitarization},
there exists an open set $W\subset\bbR^3$ such that for all
triples $(w_0,\,w_1,\,w_\infty)\in W$,
the three functions $\mu_k(w_k,\,\lambda)$ satisfy the spherical triangle
inequalities~\eqref{eq:striangle} on $\bbS^1$,
and satisfy the
strict spherical triangle inequalities on $\bbS^1\setminus\{1\}$.

By~\cite{Goldman_1998} (see also~\cite{Schmitt_Kilian_Kobayashi_Rossman_2006}),
$M_0$, $M_1$, $M_\infty$ are pointwise unitarizable
on $\bbS^1\setminus\{1\}$.
By the Unitarization Theorem~\cite{Schmitt_Kilian_Kobayashi_Rossman_2006},
the monodromy group generated by $M_0$, $M_1$, $M_\infty$ is
simultaneously $r$-unitarizable by some $c\in\LoopR$.

By \autoref{lem:monodromy-descent}, each element of the monodromy
group upstairs is a power of an element of the monodromy group downstairs.
Hence $c$ unitarizes the monodromy group upstairs.

Let $F$ be the unitary factor in the $r$-Iwasawa factorization of $C\Phi$.
Then since the monodromy group for $C\Phi$ is $r$-unitary,
then for each deck transformation of $\widetilde\Sigma_1$,
the monodromy of $C\Phi$ is the same as that of $F$.

We also have the easy closing conditions for $\Phi$, hence $F$.
Hence $\Phi_1$ induces a CMC immersion of the punctured sphere into $\bbR^3$.
\end{proof}

\begin{remark}
It is conjectured that the ends of the immersions constructed in
\autoref{thm:nnoid-symmetry} are asymptotic to half-Delaunay surfaces
(work in progress~\cite{Kilian_Rossman_Schmitt_2006}).
This is probably implied by the fact that
(1) at each end, the potential is a perturbation of the potential
for a Delaunay surface, and hence
(2) the eigenvalues of the end monodromies are equal to those of
Delaunay surfaces.
The conjecture is corroborated by the numerically computed images.
\end{remark}

\begin{remark}
A non-loop version of this proof can be used to
construct CMC $H=1$ $n$-noids in hyperbolic space.
See~\cite{Bobenko_Pavlyukevich_Springborn_2003, Bobenko_Pavlyukevich_2005} for the trinoid
and cyclic group cases.
\Note{Explicate this a little.}
\end{remark}

\begin{remark}
\label{rem:spaceform}
Immersions into the spherical and hyperbolic spaceforms
can also be constructed by a proof analogous
to that of \autoref{thm:nnoid-symmetry}.
For these constructions
a modification of the potential~\eqref{eq:potential} must be used:
the $(1-\lambda)^2$ term should be replaced by
$(\lambda-\lambda_0)(\lambda-\lambda_0^{-1})$.
For spherical space, $\lambda_0\in \bbS^1\setminus\{\pm 1\}$ 
and for hyperbolic space, $\lambda_0\in\bbR\setminus\{\pm 1\}$.
See~\cite{Schmitt_Kilian_Kobayashi_Rossman_2006} for a similar construction
in the three spaceforms carried out for trinoids.
\end{remark}

\typeout{== symmetry.tex =============================================}
\section{Symmetry}
\label{sec:symmetry}

Now the astronoid families are constructed,
it remains to compute their symmetry groups.
In this section we show that for a potential of the form~\eqref{eq:potential},
the symmetry group of the Hopf differential
induces a symmetry group of the immersion,
under certain assumptions on the monodromy group.
The proof divides into two parts:
(1) For each member $f$ of the associate family of immersions,
$G$ induces a group of isometries of Eucildean space
under which $f$ is equivariant (\autoref{sec:symmetry-assoc}).
(2) For a member of the associate family which is ``closed,'' that is,
an immersion of the universal cover which descends
to the base Riemann surface, $G$
is isomorphic to a symmetry group of the immersion
(\autoref{sec:symmetry-closed}).
Symmetries of CMC surfaces are also investigated in~\cite{Dorfmeister_Haak_1998}.

\subsection{Preliminary}
\label{sec:symmetry-preliminary}

To prove the main symmetry result of this section,
(\autoref{thm:symmetry-closed})
we require several preliminary lemmas.

\begin{lemma}
\label{lem:unitarizer-behavior}
Let $c\in\LoopR$ be a unitarizer constructed in Unitarization Theorem
of~\cite{Schmitt_Kilian_Kobayashi_Rossman_2006}.
Then there exists a finite subset $S\subset\bbS^1$ such that
$c$ extends holomorphically to $\bbS^1\setminus S$, and
has branched points or branched poles on $S$.
\Note{This lemma needs to be proved in detail.}
\end{lemma}

\begin{proof}
$c$ satisfies $c^\ast c = X$, where $X$ is the kernel of a certain
linear operator. Since $X$ has at worst branched poles, then so does $c$.
\end{proof}

\begin{lemma}
\label{lem:extend-h}
Let $h\in\LoopR$ and suppose that for some finite subset $S\subset\bbS^1$,
$h$ extends holomorphically to $\bbS^1\setminus S$,
taking values in $\matSU{2}{}$ there,
and has branch points or branched poles at each point of $S$.
Then $h$ extends holomorphically to $\bbS^1$ and is in $\LoopSL{}$.
\end{lemma}

\begin{proof}
This follows by the compactness of $\matSU{2}{}$.
\end{proof}

\begin{lemma}
\label{lem:unitary-conjugator}
Let $M_1,\,M_2\in\matSU{2}{}$ satisfy $[M_1,\,M_2]\ne 0$.
If $X\in\matSL{2}{\bbC}$ is such that
$XM_1X^{-1}\in\matSU{2}{}$ and $XM_2X^{-1}\in\matSU{2}{}$,
then $X\in\matSU{2}{}$.
\end{lemma}

\begin{proof}
Since $M_1$ is diagonalizable by an element of $\matSU{2}{}$,
we may assume without loss of generality that $M_1$ is diagonal.
Let $X=UT$ be the QR-factorization of $X$, so $U\in\matSU{2}{}$
and $T\in\matSL{2}{\bbC}$ is upper-triangular.
Then for $k=1,\,2$,
$TM_kT^{-1}\in\matSU{2}{}$, so $[T^\ast T,\,M_k]=0$.
Since $M_1$ is diagonal, then $T^\ast T$ is diagonal.
Since $[M_1,\,M_2]\ne 0$, then $M_2$ is not diagonal.
Hence $T^\ast T=\pm\id$, so $T=\pm\id$, so $X=U\in\matSU{2}{}$.
\end{proof}

\subsection{Symmetry of the associate family}
\label{sec:symmetry-assoc}

We show that the symmetry group of the Hopf differential
induces a group of symmetries on each member of the
associate family of immersions
of the universal cover of the punctured sphere.
In \autoref{sec:symmetry-closed}  we consider
the symmetry of a closed immersion, that is,
a member of the associate family which descends to the punctured sphere.

We start with a group $G$ of conformal and anticonformal
M\"obius transformations. (An ``anticonformal M\"obius
transformation'' is a M\"obius transformation post-composed with
complex conjugation.)
\Note{The terminology ``anticonformal M\"obius transformation'' is oxymoronic.
What is the standard terminology for these objects?}

We define a group $\widetilde{G}$ acting on $\widetilde{\Sigma}$
and a group homomorphism $q:\widetilde{G}\to G$ as follows:
define $\widetilde{G}$ as the set
of maps $\sigma:\widetilde{\Sigma}\to\widetilde{\Sigma}$
such that $\pi\circ\sigma = \tau\circ\pi$ for some $\tau\in G$,
and define $q(\sigma)=\tau$.

\begin{theorem}
\label{thm:symmetry-assoc}
Let $\xi$ be a potential of the form~\eqref{eq:potential}
on a punctured sphere $\Sigma$, and let $Q$ be the Hopf differential
on $\Sigma$ determined by $\xi$.
Let $G$ be a group of conformal and anticonformal M\"obius transformations
on the Riemann sphere
which is compatible with $Q$ in the following sense:
\begin{equation}
\label{eq:Qsymmetry}
\begin{split}
&\text{$\tau^\ast Q = Q$ for $\tau\in G$ conformal;}\\
&\text{$\ol{\tau^\ast Q} = Q$ for $\tau\in G$ anticonformal.}
\end{split}
\end{equation}
Let $\widetilde G$ the group which is the ``lift'' of $G$ to the universal
cover $\widetilde\Sigma$ as defined above.
Then there exists a group homomorphism $\rho:\widetilde G\to\euclid$
such that for all $\tau\in G$,
\[
\tau^\ast \left. f\right|_{\lambda=1}
 = \rho(\tau)\circ \left. f\right|_{\lambda=1}.
\]
\end{theorem}

\begin{proof}
Let $\pi:\widetilde\Sigma\to\Sigma$ be the universal cover of $\Sigma$.
Let $\Phi:\widetilde\Sigma\to\LoopR$ be a solution to the equation
$d\Phi = \Phi\xi$ with unitary monodromy which is not identically
reducible on $\bbS^1$.

We define the map $\rho:\widetilde G\to\euclid$
for the two cases of conformal and anticonformal M\"obius transformations.
Let $\tau\in\widetilde G$ and assume first that $\tau$ is conformal.
By \autoref{lem:gauge2}(i),
we have
\[
\tau^\ast \xi = \gauge{\xi}{\sgauge{z}{\tau}}.
\]
Then we have
\begin{equation}
\label{eq:sym1}
\tau^\ast\Phi = h_\tau \Phi \sgauge{z}{\tau}
\end{equation}
for some $z$-independent $h_\tau\in\LoopR$.

We now show that $h_\tau\in\LoopuSL{}$.
We have that $\Phi = c\Phi_0$ for some $c\in\LoopR$ and $\Phi_0\in\LoopSL{}$
($c$ is the unitarizer).
By \autoref{lem:unitarizer-behavior},
there exists a finite subset $S\subset\bbS^1$ such that
$c$ extends holomorphically to $\bbS^1\setminus S$,
and has branch points or branched poles on $S$.
Hence the same holds for $h$.
Since by hypothesis, the monodromy group $\calM_1$ on $\widetilde\Sigma$
is not identically reducible, there exist deck transformations
$\delta_1,\,\delta_2\in\Delta$
such that $M[\delta_1,\,\Phi]$ and $M[\delta_2,\,\Phi]$ are irreducible
on $\bbS^1$ except at finitely many points.
That is,
$[M[\delta_1,\,\Phi],\,M[\delta_2,\,\Phi]]\not\equiv 0$.
For $k=1,\,2$, define
$\epsilon_k = \tilde\tau\circ\delta_k\circ{\tilde{\tau}}^{-1}$.
By \autoref{lem:monodromy-conjugation}(ii) (with $\Psi=\Phi$ in that lemma),
we have
\[
M[\epsilon_k,\,\Phi] = \pm h M[\delta_k,\,\Phi]h^{-1},\quad k=1,\,2.
\]
Since by hypothesis, $\calM_1$ is unitary,
then for $k=1,\,2$, $h$ conjugates the two unitary
loops $M[\delta_k,\,\Phi]$ to unitary loops $\pm M[\epsilon_k,\,\Phi]$.
By \autoref{lem:unitary-conjugator}, $h$ takes values
in $\matSU{2}{}$ on $\bbS^1\setminus S$.
By \autoref{lem:extend-h}, $h$ extends to $\bbS^1$ and is in $\LoopuSL{}$.

Let $F$ be the $s$-unitary factor of $\Phi$ for $s\in(r,\,1)$.
Then~\eqref{eq:sym1} implies that $\tau^\ast F = h_\tau F$.
The Sym formula then gives
\[
\tau^\ast f = h_\tau f h_\tau^{-1} + h_\tau' h_\tau^{-1}.
\]
Define $\rho(\sigma,\,\lambda)$ by
\[
\rho(\tau,\lambda)(x) =
  h_{\sigma}x h_{\sigma}^{-1} + h_{\sigma}'h_{\sigma}^{-1}.
\]
Then by definition of $\rho(\sigma,\,\lambda)$ we have equivariance
$\tau^\ast f = \rho(\sigma,\,\lambda) \circ f$.

In a similar way we define $\rho(\tau,\,\lambda)$
for the case of anticonformal $\tau\in G$.
By \autoref{lem:gauge2}(ii),
\[
\tau^\ast \barbar{\xi} = \gauge{\xi}{\sgauge{z}{\ol\tau}},
\]
so
\[
\tau^\ast \barbar\Phi = h_\tau \Phi \sgauge{z}{\ol\tau}
\]
for some $z$-independent $h_{\sigma}\in\LoopR$.
By an argument similar to that in the conformal case above, we have
$h_\sigma\in\LoopuSL{}$.

Let $F$ be the $s$-unitary factor of $\Phi$ for $s\in(r,\,1)$.
Then~\eqref{eq:sym1} implies that $\tau^\ast \barbar{F} = h_\tau F$.
The Sym formula then gives
\[
\tau^\ast f = -\left(\barbar{h_\tau} \barbar{f} \barbar{h_\tau}^{-1} + \barbar{h_\tau}' \barbar{h_\tau}^{-1}\right).
\]
Define $\rho(\tau,\,\lambda)$ by
\[
\rho(\tau,\lambda)(x) =
-\left(\barbar{h_\tau} \barbar{x} \barbar{h_\tau}^{-1} + \barbar{h_\tau}' \barbar{h_\tau}^{-1}\right).
\]
Then by definition of $\rho(\sigma,\,\lambda)$ we have equivariance
$\tau^\ast f = \rho(\sigma,\,\lambda) \circ f$.

That the map $\rho:\widetilde G\to\euclid$ at $\lambda=1$ defined by
$\rho(\tau) = \rho(\tau,\,1)$ is a group homomorphism
follows from the equivariance:
\[
\rho(\tau_1\circ\tau_2)\circ f = f\circ\tau_1\circ\tau_2 = \rho(\tau_1)\circ f\circ \tau_2 = \rho(\tau_1)\circ\rho(\tau_2)\circ f.
\]
Since the image of $f$ spans $\matsu{2}{}$, then
$\rho(\tau_1\circ\tau_2) = \rho(\tau_1)\circ\rho(\tau_2)$.
\end{proof}

\subsection{Symmetries of closed surfaces}
\label{sec:symmetry-closed}

Equipped with the symmetry result for the associate family
(\autoref{thm:symmetry-assoc}),
we now investigate the symmetries of closed surfaces.

\begin{theorem}[Symmetry Theorem]
\label{thm:symmetry-closed}
Let $\xi$ be a potential of the form~\eqref{eq:potential}
on a punctured sphere $\Sigma$, and let $Q$ be the Hopf differential
on $\Sigma$ determined by $\xi$.
Let $G$ be a group of conformal and anticonformal M\"obius transformations
on the Riemann sphere
which is compatible with $Q$ as in~\eqref{eq:Qsymmetry}.
Suppose the symmetry group $G$ is finite.
\Note{Do we need finiteness of $G$?}
Suppose that for some $\lambda=\lambda_0$,
the CMC immersion $f_{\lambda_0}$ on the universal
cover $\widetilde\Sigma$
descends to an immersion of the punctured sphere $\Sigma$.
Then there exists a group isomorphism $\rho:G\to\euclid$
under which $f_{\lambda_0}$ is equivariant,
that is, $\tau^\ast f_{\lambda_0} = \rho(\tau)\circ f_{\lambda_0}$
for all $\tau\in G$.
\end{theorem}

\begin{proof}
We first construct the group homomorphism $\rho:G\to\euclid$.
Consider the group $\widetilde G$ defined in \autoref{sec:symmetry-assoc},
and let $\tilde\rho:\widetilde G\to\euclid$ be the
group homomorphism constructed in \autoref{thm:symmetry-assoc}.

\setlength{\unitlength}{60pt}
\begin{figure}[ht]
\begin{picture}(1.0,1.2)(0,0)
\put(0,1){$\widetilde{G}$}
\put(0.15,1.05){\vector(1,0){0.85}}
\put(0.5,1.12){$\tilde\rho$}
\put(0,0){$G$}
\put(0.05,0.95){\vector(0,-1){0.8}}
\put(-0.1,0.5){$q$}
\put(0.15,0.05){\vector(1,1){0.9}}
\put(0.6,0.4){$\rho$}
\put(1,1){$T$}
\end{picture}
\end{figure}

Let $\Delta$ be the group of deck transformations of $\widetilde{\Sigma}$.
Note that $\Delta$ is a subgroup of $\widetilde{G}$,
We now show that $\Delta\subset\ker\tilde\rho$.
The form of the potential $\xi$ implies the closing
conditions~\eqref{eq:closing}.
Let $\delta\in\Delta$.
Then $\delta^\ast\Phi = h_\delta\Phi$,
where $h_\delta$ is the monodromy of $\Phi$ corresponding to $\delta$.
Then $h_\delta(1)=\id$ and $h_\delta'(1)=0$, so
$\rho(\delta)=\id$.

This shows that $\Delta\subset\ker\tilde\rho$. It follows that there a
group homomorphism $\rho:G\to\euclid$
so that $\tilde\rho = q \circ \rho$,
that is, the above diagram commutes.

We now show that $\rho$ is a group isomorphism.
Let $\tau\in G$ and define $h_\tau$ by $\tau^\ast\Phi = h_\tau \Phi g$,
$g =\sgauge{z}{\tau}$.
By \autoref{lem:monodromy-descent}(i),
$h$ is the monodromy downstairs.
The eigenvalues of the monodromy downstairs
is $\rho$.
Then $\rho(1) = \half-\tfrac{1}{2n}$,
where $n\in\bbN_{>1}$ is the order of $u$ at a preimage of a branch point.
Hence $\rho(1)\not\in\half\bbZ$.
Hence $M(1)\not\in\{\pm\id\}$,
so $h(1)\not\in\{\pm\id\}$.
Hence $\rho(\tau)$ is not the trivial element of $\euclid$.
This shows that $\ker\rho$ is trivial, so $\rho$ is an isomorphism.
\end{proof}

\begin{remark}
The symmetry \autoref{thm:symmetry-assoc} and \autoref{thm:symmetry-closed}
apply when the immersion is constructed by $r$-dressing,
as long as the dressing can be extended holomorphically to the annulus
$\{\lambda\in\bbC\suchthat r<\abs{\lambda}<1\}$.
If the dressing cannot be so extended, 
then it is not known whether $h$ in the proof extends
holomorphically to $\bbS^1$,
a requirement for the completion of the proof.
Thus the theorem does not apply
to $n$-noids constructed by dressing by simple factors;
images of theses surfaces~\cite{Kilian_Schmitt_Sterling_2004} indicate that dressing breaks
symmetry.
\end{remark}

\begin{remark}
The symmetries of $n$-noid immersions into
the spherical and hyperbolic spaceforms
(\autoref{rem:spaceform})
can be computed by an argument analogous to that of
\autoref{thm:symmetry-closed}.
\end{remark}

\section{Symmetries of astronoids}

We conclude by computing the symmetries of the $n$-noids
constructed in \autoref{thm:nnoid}.
For each choice of finite subgroup $G\subset\matPSL{2}{\bbC}$,
let $H$ be the corresponding subgroup of $\matSO{3}{}$.
Then the symmetry group of an $n$-noid constructed from $G$
is generated by the elements of $H$ together with a
reflection in a plane.

\typeout{== irreducible.tex =============================================}
\subsection{Irreducibility}

As the irreducibility of the monodromy representation is a requirement
for the application of the Symmetry \autoref{thm:symmetry-closed},
we show that the unitary monodromy group of the $n$-noids constructed
in \autoref{thm:nnoid} is irreducible.
The proof divides into cases according to how many of the three end weights
are zero.
For the case of two or three non-zero weights,
an argument can be made using the fact
that the monodromies are powers of the irreducible ``trinoid''
monodromies group downstairs.
For the case of one non-zero weight, a special argument must be used.
If all three weights are zero, the $n$-noid degenerates to a sphere.

We start with a simple fact about the irreducibility of powers
of matrices.

\begin{lemma}
\label{lem:irreducible-powers}
Let $A,\,B\in\mattwo{\bbC}$ be irreducible. Suppose for
some $r,\,s\in\bbN_{>0}$,
$A^r\not\in\bbC\id$ and $B^s\not\in\bbC\id$.
Then $A^r$ and $B^s$ are irreducible.
\end{lemma}

\begin{proof}
If $A^r$ and $B^s$ are reducible, they have a common eigenline $v$.
Since $v$ is an eigenline of $A^r$, then $v$ is an eigenline $A$.
Similarly $v$ is an eigenline of $B$. Thus $A$ and $B$ have a common
eigenline, and are hence reducible.
\end{proof}

\begin{theorem}
\label{thm:irreducible}
Let $u$, $\xi$ and $\eta$ be as in \autoref{not:nnoid}.
Assume that at least one of the three weights $w_0,\,w_1,\,w_2$
definining $Q$ is non-zero.
Let $\Phi$ be the meromorphic solution
to the equation $d\Phi=\Phi\xi$ as in \autoref{thm:nnoid},
with initial condition chosen so that $\Phi$ has unitary monodromy group
$\calM$.
Then $\calM$ is irreducible on $\bbS^1$ except possibly on a finite
subset.
\end{theorem}

\begin{proof}
\textit{Case 1: two or three non-zero weights:}
Let $M_1$ and $M_2$ be monodromies downstairs
around ends with non-zero weights.
By \autoref{lem:trinoid}, the monodromy group downstairs
is irreducible on $\bbS^1$ except at a finite subset.
By \autoref{lem:irreducible-powers}, the same holds for the
monodromy group upstairs.

\textit{Case 2: one non-zero weight:}
Arrange so that the end with non-zero weight is the orbit of
$\tau:z\mapsto \alpha z$, where $\alpha$ is a root of unity.
Then $\tau^\ast\xi = \gauge{\xi}{g}$,
where $g=\diag(\alpha^{1/2},\,\alpha^{-1/2})$.
Then $0$ is a fixed point of $\tau$, and is not a pole of $\xi$.

We now use an argument found in~\cite{Rossman_Schmitt_2006},
there applied to the cyclic case, and here applied more generally.
We have $M_k = g^kM_0g^{-k}$. We also have $M_0\cdots M_{n-1}=\id$.
Hence $(M_0g)^n=-\id$. It follows that the eigenvalues of $M_0 g$ are
$n$'th roots of $-1$, and are hence constant.
Using $M_0(1)=\id$, we get that the eigenvalues of $M_0 g$ are
$\alpha^{\pm 1/2}$.

We now apply \autoref{lem:trinoid} to the triple
$(M_0,\,g,\,(M_0g)^{-1})$. Let the eigenvalues of $M_0=\exp(\pm 2\pi i\nu)$,
$\nu\in[-1/2,\,1/2]$.
By \autoref{lem:trinoid},
$(M_0,\,g,\,(M_0g)^{-1})$ are irreducible on $\bbS^1$ if $0<\abs{\nu}<1/n$.
Since $\nu$ is not constant on $\bbS^1$,
the inequality holds on $\bbS^1$ except possibly at a finite subset.
\end{proof}

\typeout{== main.tex =============================================}
\subsection{\texorpdfstring{Symmetric $n$-noids}{Symmetric n-noids}}
\label{sec:main}

The final theorem below computes the symmetry groups of the $n$-noids
constructed in \autoref{thm:nnoid}.
It refers to the following finite subgroups of $\matO{3}{}$,
all of which contain orientation-reversing elements:
\begin{itemize}
\item
$n$-fold pyramidal symmetry is the symmetry of a right pyramid over
a regular $n$-gon, with group order $2n$.
\item
$n$-fold prismatic symmetry is the symmetry of a right prism
(or double-pyramid) over a regular $n$-gon, with group order $4n$.
\item
Full tetrahedron, full octahedron and full icosahedron symmetries are
the full symmetry groups of the respective regular polyhedra,
with respective group orders $24$, $48$ and $120$.
\end{itemize}

\begin{theorem}
\label{thm:nnoid-symmetry}
Let $G\subset\matPSL{2}{\bbC}$ be a finite subgroup,
and let $f$ be a CMC immersion of the punctured Riemann sphere
constructed in \autoref{thm:nnoid}.
Then $f$ has the following symmetries:
for $G=\bbZ$, $f$ has $n$-fold pyramidal symmetry;
for $G=D_n$, $f$ has $n$-fold prismatic symmetry,
and for $G=A_4$, $S_4$ and $A_4$, $f$ has respectively
full tetrahedral, full octahedral or full icosohedral symmetry.
\end{theorem}

\begin{proof}
Let $H$ be the subgroup of the conformal and anticonformal M\"obieus
transforms generated by $G$ together with the map $z\mapsto\ol{z}$.
The proof uses the Symmetry \autoref{thm:symmetry-closed},
with $H$ playing the role of $G$ in that theorem.

By the choice of $Q$, for all $\tau\in G$ we have $\tau^\ast Q = Q$.
Let $c=\conj$.
We now show that $\ol{c^\ast Q}=Q$.
We have
$Q=u^\ast q$ and $\ol{c^\ast q} = q$.
Moreover, $u$ was chosen so that
$u\circ c = c\circ u$, so $c^\ast u^\ast = u^\ast c^\ast$.
Then
\[
\ol{c^\ast Q} = \ol{c^\ast u^\ast q}
 = \ol{ u^\ast c^\ast q}
 = u^\ast q = Q.
\]
Hence $Q$ has the symmetries with respect to $H$ as required by
\autoref{thm:symmetry-closed}.

By \autoref{lem:irreducible-powers}, the same holds for the
monodromy group upstairs.

By \autoref{thm:irreducible}, the monodromy group of $\Phi$ is
irreducible.
\Note{Explicate this proof.}
\end{proof}



\bibliographystyle{amsplain}
\bibliography{symmetric-nnoid}

\end{document}